\documentclass[12pt]{amsart}
\usepackage{amsmath, amsfonts, amssymb, amsthm,hyperref,mathtools,array}
\hypersetup{hypertex=true,
	colorlinks=true,
	linkcolor=blue,
	anchorcolor=blue,
	citecolor=blue}
\usepackage{bm}
\allowdisplaybreaks[4]
\def\({\bg(}
\def\){\bg)}

\def\ord{{\rm ord}}

\def\gen{{\rm gen}}
\def\spn{{\rm spn}}

\def\Z{\mathbb Z}

\def\N{\mathbb N}
\def\Q{\mathbb Q}

\def\rank{{\rm rank}}

\def\Norm{{\rm Norm}}
\def\char{{\rm char}}

\def\v{{\bm v}}

\def\x{{\bm x}}

\def\e{{\bm e}}
\def\0{{\bm 0}}
\def\s{{\mathfrak{s}}}

\def\Hom{{\rm Hom}}
\def\char{{\rm char}}
\def\pmod #1{\ ({\rm{mod}}\ #1)}
\def\mod #1{\ {\rm mod}\ #1}
\def\Ack{\medskip\noindent {\bf Acknowledgments}}

\theoremstyle{plain}
\newtheorem{theorem}{Theorem}[section]
\newtheorem{lemma}{Lemma}
\newtheorem{corollary}{Corollary}
\newtheorem{conjecture}{Conjecture}

\theoremstyle{definition}

\theoremstyle{remark}
\newtheorem{remark}{Remark}

\makeatletter
\@namedef{subjclassname@2020}{%
	\textup{2020} Mathematics Subject Classification}
\makeatother
\vspace{4mm}

\begin{document}
	\medskip
	
	\title[On restricted sums of four squares and Zhi-Wei Sun's $x+24y$ conjecture]
	{On restricted sums of four squares and Zhi-Wei Sun's $x+24y$ conjecture}
	\author[H.-L. Wu and Y.-F. She]{Hai-Liang Wu and Yue-Feng She*}
	
	\address {(Hai-Liang Wu) School of Science, Nanjing University of Posts and Telecommunications, Nanjing 210023, People's Republic of China}
	\email{\tt whl.math@smail.nju.edu.cn}
	
	\address {(Yue-Feng She) Department of Applied Mathematics, Nanjing Forestry University, Nanjing 210037, People's Republic of China}
	\email{\tt she.math@njfu.edu.cn}

	\keywords{restricted sums of four squares, quadratic forms. 
		\newline \indent This research was supported by the Natural Science Foundation of China (Grant No. 12101321).
	    \newline \indent *Corresponding author.}
	
	\begin{abstract}
	In this paper, by using the arithmetic theory of ternary quadratic forms, we study some refinements on Lagrange's four-square theorem. For example, given positive integers $a,b$ satisfying some algebraic conditions and a positive integer $C\ge3$, we will show that for any sufficiently large integer $n$ with $\ord_2(n)\le C$, there exist non-negative integers $x,y,z,w$ such that 
	$$\begin{cases}
		x^2+y^2+z^2+w^2=n,\\
		ax+by\in\mathcal{S},
	\end{cases}$$
    where $\mathcal{S}$ is the set of all squares over $\mathbb{Z}$. In particular, we obtain some progress on Zhi-Wei Sun's $x+24y$ conjecture.
	\end{abstract}
	\maketitle

	\section{Introduction}
	\setcounter{lemma}{0}
	\setcounter{theorem}{0}
	\setcounter{equation}{0}
	\setcounter{conjecture}{0}
	\setcounter{remark}{0}
	\setcounter{corollary}{0}

    \subsection{Background and Motivation}

     In 1770, Lagrange established the celebrated four-square theorem (cf. \cite[pp. 5--7]{Nathanson}), which states that each number $n\in\mathbb{N}=\{0,1,2,\cdots\}$ can be written as 
     $$n=x^2+y^2+z^2+w^2$$
     with $x,y,z,w\in\mathbb{Z}$, where $\mathbb{Z}$ denotes the ring of integers. Along this line, in 1917, Ramanujan \cite{Ramanujan} proved that there are at most $55$ possible quadruples
     $(a,b,c,d)$ of positive integers with $a\le b\le c\le d$ such that the quadratic form $ax^2+by^2+cz^2+dw^2$ is {\it universal} over $\mathbb{Z}$, i.e., 
     $$\left\{ax^2+by^2+cz^2+dw^2:\ x,y,z,w\in\mathbb{Z}\right\}=\mathbb{N}.$$
     Ten years later, Dickson \cite{Dickson} confirmed that $54$ quadratic forms on Ramanujan's list are indeed universal over $\mathbb{Z}$, and pointed out that for the remaining one quadruples $(1,2,5,5)$, we have 
     $$\left\{x^2+2y^2+5z^2+5w^2:\ x,y,z,w\in\mathbb{Z}\right\}=\mathbb{N}\setminus\{15\}.$$

     Cauchy first considered a refinement of Lagrange's four-square theorem and obtained the following result, which was later known as Cauchy's lemma (cf. \cite[p. 31]{Nathanson}). Let $n$ and $m$ be positive odd integers with $m^2<4n$ and $3n<m^2+2m+4$. Then there exist $x,y,z,w\in\mathbb{N}$ such that 
     $$\begin{cases}
     		n=x^2+y^2+z^2+w^2,\\
     		m=x+y+z+w.
     \end{cases}$$
     Along this line, for positive integers $m,n,t$, Kloosterman \cite{Kloosterman} considered 
     \begin{equation*}
     r_t(n,m):=\#\left\{(x_1,x_2,\ldots,x_t)\in \mathbb{Z}^t:\ 
     	x_1^2+x_2^2+\cdots+x_t^2=n\ \text{and}\ x_1+x_2+\cdots+x_t=m\right\},
     \end{equation*}
     and obtained an asymptotic formula via circle method, 
     where $\# S$ denotes the cardinality of a finite set $S$. Later Lomadze \cite{Lomadze} and Van der Blij \cite{Blij} also obtained some developments on this topic.

    Inspired by the above results, in 2017, Z.-W. Sun \cite{Sun17} initiated the study of the following refinement of Lagrange's four-square theorem. Let $P(X,Y,Z,W)\in\mathbb{Z}[X,Y,Z,W]$. Then by Sun's definition, the polynomial $P(X,Y,Z,W)$ is said to be {\it suitable}, if for any $n\in\mathbb{N}$, there are $x,y,z,w\in\mathbb{N}$ such that 
    $$\begin{cases}
    	n=x^2+y^2+z^2+w^2,\\
    	P(x,y,z,w)\in\mathcal{S},
    \end{cases}$$
    where 
    $$\mathcal{S}=\{x^2:\ x\in\mathbb{Z}\}$$ 
    is the set of all squares over $\mathbb{Z}$. 
    
    The most well-known problem in this topic is Sun's 1-3-5 conjecture (see \cite[Conjecture 4.3(i)]{Sun17}). This challenging conjecture states that for any integer $n\in\mathbb{N}$, there exist $x,y,z,w\in\mathbb{N}$ such that 
    $$\begin{cases}
    	n=x^2+y^2+z^2+w^2,\\
    	x+3y+5z\in\mathcal{S},
    \end{cases}$$
     Later Machiavelo and his collaborators \cite{MT,MT2} confirmed this conjecture by using quaternion arithmetic in the ring of Lipschitz integers. Furthermore, for $n\in\mathbb{N}$, $m\in\mathbb{Z}$ and $\v=(a,b,c,d)\in\mathbb{Z}^4$, we use $r_{\v}(n,m)$ or $r_{(a\text{-}b\text{-}c\text{-}d)}(n,m)$ to denote the number of solutions $(x,y,z,w)\in\Z^4$ for the system
    $$\begin{cases}
    	n=x^2+y^2+z^2+w^2,\\
    	m=ax+by+cz+dw.
    \end{cases}$$  
     In 2024, by using the theory of Siegel modular forms, D. Li and H. Zhou \cite{LiZhou} gave an explict formula for 
      $$r_{(1\text{-}3\text{-}5\text{-}0)}(n,m)+r_{(1\text{-}3\text{-}3\text{-}4)}(n,m),$$ 
      which involves the Hurwitz class numbers. Later, for any $t\in\N$, J. Li and H. Wang \cite{LiWang} obtained a formula for 
      $$\sum_{a^2+b^2+c^2+d^2=t}r_{\v}(n,m).$$ 
      However, it remains difficult to compute the explicit value of $r_{\v}(n,m)$ separately for general $\v\in\mathbb{Z}^4$.

    On the other hand when $P(X,Y,Z,W)=aX+bY$, Sun \cite[Conjecture 4.1]{Sun17} posed the following conjecture.
     
     \begin{conjecture}[Z.-W. Sun]
     	Let $a,b\in\Z^+$ with $\gcd(a,b)$ square-free and $a\le b$. Then $aX+bY$ is  suitable if and only if 
     	$$(a,b)=(1,2),(1,3),(1,24).$$ 
     \end{conjecture}

      For this conjecture, some partial results were obtained. 
 
     \begin{itemize}
     	\item Y.-C. Sun and Z.-W. Sun \cite[Theorem 1.1(ii)]{Sun-aa} showed that $X+2Y$ is suitable.
     	
     	\item Y.-F. She and H.-L. Wu \cite[Theorem 1.4]{She-Wu} proved that $X+3Y$ is  suitable.
     \end{itemize}
     In addition, as mentioned by Sun, confirming the suitability of $X+24Y$ might be very difficult.
     
     Motivated by the above results, in this paper, we will concentrate on the case when the polynomial $P(X,Y,Z,W)$ is a binary homogeneous integral polynomial of degree one.

     \subsection{Main Results}  For any positive integer $l\ge2$, let 
     $$\left(\mathbb{Z}/l\mathbb{Z}\right)^{\times}=\left\{x \mod l\mathbb{Z}:\ x\in\mathbb{Z}\ \text{and}\ \gcd(x,l)=1\right\}$$
     be the multiplicative group of all units in the residue ring $\mathbb{Z}/l\mathbb{Z}$. Now we state our main result. 

     \begin{theorem}\label{Thm. A}
     	Let $a,b\in\mathbb{Z}^+$ with $\gcd(a,b)=1$ such that $(\mathbb{Z}/(a^2+b^2)\mathbb{Z})^{\times}$ is a cyclic group. 
     	Let $C\ge 3$ be a positive integer. Then for any sufficiently large integer $n$ with $\ord_2(n)\le C$, there exist $x,y,z,w\in\mathbb{N}$ such that 
     	$$\begin{cases}
     		x^2+y^2+z^2+w^2=n,\\
     		ax+by\in\mathcal{S}.
     	\end{cases}$$
     \end{theorem}

      \begin{remark}
      	(i) Clearly $(\mathbb{Z}/(a^2+b^2)\mathbb{Z})^{\times}$ is a cyclic group if and only if 
      	\begin{equation}\label{Eq. A in Thm. A}
      		a^2+b^2\in\left\{p^r:\ r\in\mathbb{Z}^+\ \text{and}\ p\ \text{is an odd  prime}\right\}\cup\left\{2p^k:\ k\in\mathbb{N}\ \text{and}\ p\ \text{is an odd prime}\right\}.
      	\end{equation}

      	(ii) The condition $\ord_2(n)\le C$ can not be removed. For example, let $(a,b)=(3,10)$, which satisfies the condition (\ref{Eq. A in Thm. A}). It is known that each positive integer of the form $4^{2r+1}\cdot 6$ has a unique partition into four integral squares (see \cite[p. 86]{Gro}), that is, 
      	$$4^{2r+1}\cdot 6=(4^{r+1})^2+(2\cdot 4^r)^2+(2\cdot 4^r)^2+0^2.$$
      	However, it is easy to verify that for any positive integer of the form $4^{2r+1}\cdot 6$, there are no $x,y,z,w\in\mathbb{N}$ such that 
      	$$\begin{cases}
      		x^2+y^2+z^2+w^2=4^{2r+1}\cdot 6,\\
      		3x+10y\in\mathcal{S}.
      	\end{cases}$$
      \end{remark}

       As $(a,b)=(1,24)$ satisfies the condition (\ref{Eq. A in Thm. A}), we directly obtain the following corollary, which makes some progress on Sun's $x+24y$ conjecture.
       
       \begin{corollary}\label{Cor. x+24y conjecture}
       	Let $C\ge 3$ be a positive integer. Then for any sufficiently large integer $n$ with $\ord_2(n)\le C$, there exist $x,y,z,w\in\mathbb{N}$ such that 
       	$$\begin{cases}
       		x^2+y^2+z^2+w^2=n,\\
       		x+24y\in\mathcal{S}.
       	\end{cases}$$
       \end{corollary}

      \subsection{Outline of this paper} In Section 2, we will introduce some basic notations and necessary lemmas concerning lattices on quadratic spaces. The proof of our main result will be given in Section 3.

     \section{Preliminaries concerning lattices on quadratic spaces}
      \setcounter{lemma}{0}
      \setcounter{theorem}{0}
      \setcounter{equation}{0}
      \setcounter{conjecture}{0}
      \setcounter{remark}{0}
      \setcounter{corollary}{0}
     
     \subsection{Notation} 
     
     We first adopt standard notations appeared in \cite{OTO}. Let $R$ be a principal ideal domain and let $\mathbb{F}$ be the fractional field of $R$ with $\char({\mathbb{F}})\neq 2$, where $\char({\mathbb{F}})$ is the characteristic of $\mathbb{F}$. Let $(U,B,Q)$ be a non-degenerated quadratic space over $\mathbb{F}$ with $\dim_{\mathbb{F}}U=m$, where $B$ is a symmetric bilinear form on $U$ and $Q$ is the corresponding quadratic map.	The group of all rotations of $U$ is defined by 
     $$O^+(U):=\left\{\sigma\in\Hom_{\mathbb{F}}(U,U):\ \det(\sigma)=1\ \text{and}\ Q(\sigma(\x))=Q(\x)\ \text{for any}\ \x\in U\right\}.$$
     Also, let 
     $$\theta_{\mathbb{F}}:\ O^+(U)\longrightarrow \mathbb{F^{\times}}/\mathbb{F}^{\times2}$$ be the spinor norm map, where 
     $$\mathbb{F}^{\times }=\mathbb{F}\setminus\{0\},$$
     and
     $$\mathbb{F}^{\times 2}=\left\{x^2:\ x\in\mathbb{F}\setminus\{0\}\right\}.$$
     The kernel of $\theta_{\mathbb{F}}$ is denoted by $O'(U)$. Let 
     $$M=R\e_1+R\e_2+\cdots+R\e_m$$
     be a non-degenerated lattice on $U$, where $\e_1,\e_2,\cdots, \e_m$ is an $R$-basis of $M$. Then the gram matrix $A$ of $M$ with respect to $\e_1,\e_2,\cdots, \e_m$ is defined by 
     $$A=\left(B(\e_i,\e_j)\right)_{1\le i,j\le m}.$$
     Conversely, for any $m\times m$ symmetric matrix $A$ over $\mathbb{F}$, we write $M\cong A$, if there is an $R$-basis of $M$ such that $A$ is the gram matrix of $M$ with respect to this basis. For any symmetric matrices $A_1,A_2,\cdots,A_s$ over $\mathbb{F}$, we simply write 
     \begin{equation}\label{Eq. M cong As}
     	M\cong A_1\perp A_2\perp\cdots\perp A_s
     \end{equation}
     if 
     $$M\cong\left(\begin{array}{cccc}
     	A_1      & 0      & \cdots &  0\\
     	0        & A_2    & \cdots &  0\\
     	\vdots   & \vdots & \ddots &  \vdots\\
     	0        & 0      & \cdots &  A_s
     \end{array}\right).$$
     In particular, if $A_i=(a_i)$ is a $1\times 1$ matrix for $i=1,2,\cdots,s$, then (\ref{Eq. M cong As}) is abbreviated as 
     $$M\cong \langle a_1 \rangle \perp \langle a_2 \rangle \cdots \perp \langle a_s \rangle.$$
     
     Throughout the remaining part of this paper, $\mathbb{Q}$ denotes the field of rational numbers. Let $\Omega_{\mathbb{Q}}$ be the set of all places of $\mathbb{Q}$. By the Ostrowski theorem, each non-archimedean place arises from the $p$-adic valuation $|\cdot|_p$, where $p$ is a rational prime. Also, the unique archimedean place of $\mathbb{Q}$, denoted by $\infty$, arises from the ordinary absolute value $|\cdot|_{\infty}$ on $\mathbb{Q}$. Hence, for simplicity, we write   
     $$\Omega_{\mathbb{Q}}=\left\{p:\ p\ \text{is a rational prime}\right\}\bigcup\left\{\infty\right\}.$$
     
     For any $p\in\Omega_{\mathbb{Q}}\setminus\{\infty\}$, let $\mathbb{Q}_p$, $\mathbb{Z}_p$ and $\mathbb{Z}_p^{\times}$ denote the $p$-adic number field, the ring of $p$-adic integers and the multiplicative group of all $p$-adic units over $\mathbb{Z}_p$ respectively. Also, for the archimedean place $\infty$, we let $\mathbb{Q}_{\infty}=\mathbb{R}$ be the field of real numbers.

     Now consider a non-degenerated quadratic space $(V,B,Q)$ over $\mathbb{Q}$ with $\dim_{\mathbb{Q}}(V)=3$. Given a non-degenerated lattice $L$ on $V$ with $\rank(L)=3$, for any $p\in\Omega_{\mathbb{Q}}\setminus\{\infty\}$, let 
     $$V_{p}:=V\otimes_{\mathbb{Q}}\mathbb{Q}_p,$$
     and
     $$L_p:=L\otimes_{\mathbb{Z}}\mathbb{Z}_p$$
     be the $p$-adic completions of $V$ and $L$ respectively. For $\infty\in\Omega_{\mathbb{Q}}$, let 
     $$V_{\infty}=L_{\infty}=V\otimes_{\mathbb{Q}}\mathbb{Q}_{\infty}=V\otimes_{\mathbb{Q}}\mathbb{R}.$$
     The adelization of $O^+(V)$ is defined by 
     $$O_{\mathbb{A}}^+(V):=\left\{\sigma=(\sigma_v)\in\prod_{v\in\Omega_{\mathbb{Q}}}O^+(V_v):\ \sigma_v(L_v)=L_v\ \text{for almost all}\ v\in\Omega_{\mathbb{Q}}\right\}.$$
     Clearly $O^+(V)$ can be viewed as a subgroup of $O_{\mathbb{A}}^+(V)$. 
     Also, let 
     $$O_{\mathbb{A}}'(V):=\left\{\sigma=(\sigma_v)\in O_{\mathbb{A}}^+(V):\ \sigma_v\in O'(V_v)\ \text{for any}\ v\in\Omega_{\mathbb{Q}}\right\},$$
     and 
     $$O_{\mathbb{A}}^+(L):=\left\{\sigma=(\sigma_v)\in O_{\mathbb{A}}^+(V):\ \sigma_v(L_v)=L_v\ \text{for any}\ v\in\Omega_{\mathbb{Q}}\right\}.$$
     Then there is a group action of $O_{\mathbb{A}}^+(V)$ on the set of all lattices in the genus of $L$. It is also known that 
     $$\gen(L)=O_{\mathbb{A}}^+(V)L,$$
     and
     $$\spn(L)=\spn^+(L)=O^+(V)O_{\mathbb{A}}'(V)O_{\mathbb{A}}^+(L)L,$$
     where $\gen(L)$ denotes the set of all lattices in the genus of $L$, $\spn(L)$ is the set of all lattices in the spinor genus of $L$ and $\spn^+(L)$ is the set of all lattices in the proper spinor genus of $L$(since $\rank(L)=3$, we clearly have $\spn(L)=\spn^+(L)$). 
     
     For any rational number $a\in\mathbb{Q}$, we say that $a$ is represented by $\gen(L)$, denoted by $a\in Q(\gen(L))$, if there is a lattice $M\in \gen(L)$ such that 
     $$a\in Q(M):=\left\{Q(\x):\ \x\in M\right\}.$$
     It is known that 
     $$a\in Q(\gen(L))\Leftrightarrow a\in Q(L_v):=\left\{Q(\x):\ \x\in L_v\right\}\ \text{for any}\ v\in\Omega_{\mathbb{Q}}.$$
     Also, we say that $a$ is represented by $\spn(L)$, denoted by $a\in Q(\spn(L))$, if there exists a lattice $M\in\spn(L)$ such that $a\in Q(M)$. 

    \subsection{Some necessary lemmas}

     Let $p\in\Omega_{\mathbb{Q}}\setminus\{\infty\}$ be a place and $(U,B,Q)$ be a non-degenerated ternary quadratic space over $\mathbb{Q}_p$. We say that $U$ is isotropic if there exists an $\x\in U\setminus\{\0\}$ such that $Q(\x)=0$. Otherwise, $U$ is said to be anisotropic.
     
     For any non-degenerated ternary lattice 
     $$M=\mathbb{Z}_p\e_1+\mathbb{Z}_p\e_2+\mathbb{Z}_p\e_3\subseteq U,$$
     the scale of $M$, denoted by $\s(M)$, is a fractional ideal of $\mathbb{Z}_p$ generated by the set 
     $$\left\{B(\e_i,\e_j):\ 1\le i,j\le 3\right\}.$$
     Then the lattice $M$ is said to be a unimodular lattice over $\mathbb{Z}_p$ if $\s(M)=\mathbb{Z}_p$ and $d(M)\in\mathbb{Z}_p^{\times}\mod \mathbb{Z}_p^{\times 2}$, where $d(M)$ is the discriminant of $M$ and $\mathbb{Z}_p^{\times 2}=\left\{x^2:\ x\in\mathbb{Z}_p^{\times}\right\}$.

    For unimodular lattices over $\mathbb{Z}_p$, we have the following results (see \cite[92:1 and 93:11]{OTO}). 

     \begin{lemma}\label{Lem. unimodular lattices over Zp}
     	Let $p\in\Omega_{\mathbb{Q}}\setminus\{\infty\}$ be a place and $(U,B,Q)$ be a non-degenerated ternary quadratic space over $\mathbb{Q}_p$. Let 
     	$$M\cong \langle\varepsilon_1\rangle \perp \langle\varepsilon_2\rangle \perp \langle\varepsilon_3\rangle\subseteq U$$
     	 be a ternary unimodular lattice over $\mathbb{Z}_p$ with $\varepsilon_1,\varepsilon_2,\varepsilon_3\in\mathbb{Z}_p^{\times}$. 
     	
     	{\rm (i)} Suppose $p\neq 2$. Then 
     	$$M\cong\langle1\rangle \perp \langle1\rangle \perp\langle d(M)\rangle.$$
     	Moreover, if $U$ is isotropic, then 
     	$$M\cong \left(\begin{array}{cc}
     		0 & 1\\
     		1 & 0
     	\end{array}\right)\perp\langle-d(M)\rangle.$$
     
     {\rm (ii)} Suppose $p=2$. If $U$ is isotropic, then 
     $$M\cong \left(\begin{array}{cc}
     	0 & 1\\
     	1 & 0
     \end{array}\right)\perp\langle-d(M)\rangle.$$
     If $U$ is anisotropic, then 
     $$M\cong \left(\begin{array}{cc}
     	2 & 1\\
     	1 & 2
     \end{array}\right)\perp\langle-3d(M)\rangle.$$
     Moreover, 
     \begin{align*}
     	 Q\left(\left(\begin{array}{cc}
     		2 & 1\\
     		1 & 2
     	\end{array}\right)\right)
     &=\left\{2x^2+2xy+2y^2:\ x,y\in\mathbb{Z}_2\right\}\\
     &=\left\{x\in\mathbb{Z}_2:\ \ord_2(x)\equiv 1\pmod {2}\right\}\bigcup\left\{0\right\}.
     \end{align*}
     \end{lemma}
	
	Let $v\in\Omega_{\mathbb{Q}}$. Then the Hilbert symbol 
	$$(\cdot,\cdot)_v:\ \mathbb{Q}_v^{\times}/\mathbb{Q}_v^{\times 2}\times \mathbb{Q}_v^{\times}/\mathbb{Q}_v^{\times 2} \longrightarrow \left\{1,-1\right\}$$
	is defined by 
	$$(a,b)_v=\begin{cases}
		1  & \mbox{if}\ a \in \Norm_{\mathbb{Q}_v(\sqrt{b})/\mathbb{Q}_v}(\mathbb{Q}_v(\sqrt{b})^{\times}),\\
		-1 &  \mbox{otherwise,}
	\end{cases}$$
	where 
	$$\Norm_{\mathbb{Q}_v(\sqrt{b})/\mathbb{Q}_v}(\mathbb{Q}_v(\sqrt{b})^{\times})=\left\{x^2-by^2:\ x,y\in\mathbb{Q}_v\right\}\setminus\{0\}.$$
	
	We also need the following result (see \cite[58:6]{OTO}).
	
	\begin{lemma}\label{Lem. isotropic}
		Let $v\in\Omega_{\mathbb{Q}}$ and $(U,B,Q)$ be a non-degenerated ternary quadratic space over $\mathbb{Q}_v$. Let 
		$$M\cong \langle a_1\rangle \perp \langle a_2\rangle \perp \langle a_3\rangle\subseteq U$$
		be a non-degenerated lattice on $U$ with $a_1,a_2,a_3\in\mathbb{Q}_v$. Then 
		$$U\ \text{is isotropic}\ \Leftrightarrow S_v(U)=(-1,-1)_v,$$
		where 
		$$S_v(U)=\prod_{1\le i\le j\le 3}(a_i,a_j)_v$$
		is the Hasse invariant of $U$.  
	\end{lemma}
	
	 The next lemma determines the number of spinor genera in $\gen(L)$ (see \cite[102:7 and 102:9]{OTO}).
	 
	 \begin{lemma}\label{Lem. number of spinor genera}
	 	Let $(V,B,Q)$ be a non-degenerated quadratic space over $\mathbb{Q}$ with $\dim_{\mathbb{Q}}(V)=3$, and let $L$ be a non-degenerated ternary lattice on $V$. Then 
	 	\begin{align*}
	 		\#\left\{\spn(M):\ M\in\gen(L)\right\}
	 		&=\left[O_{\mathbb{A}}^+(V):\ O^+(V)O_{\mathbb{A}}'(V)O_{\mathbb{A}}^+(L)\right]\\
	 		&=\left[J_{\mathbb{Q}}:\ \mathbb{Q}^{\times}J_L\right],
	 	\end{align*}
	 	where 
	 	$$J_{\mathbb{Q}}=\left\{(i_v)_v\in\prod_{v\in\Omega_{\mathbb{Q}}}\mathbb{Q}_v^{\times}:\ |i_v|_v=1\ \text{for almost all}\ v\in\Omega_{\mathbb{Q}}\setminus\{\infty\}\right\}$$
	 	is the group of all id$\grave{e}$les,  
	 	and 
	 	$$J_L=\left\{x=(x_v)\in J_{\mathbb{Q}}:\ x_v\in\theta_{\mathbb{Q}_v}(O^+(L_v))\ \text{for any}\ v\in\Omega_{\mathbb{Q}}\right\}.$$
	 	In particular, if $\mathbb{Z}_p^{\times}\subseteq \theta_{\mathbb{Q}_p}(O^+(L_p))$ for any $p\in\Omega_{\mathbb{Q}}\setminus\{\infty\}$, then $\gen(L)=\spn(L)$. 
	 \end{lemma}
	
       We next state some necessary results to compute spinor norms for lattices over $\mathbb{Z}_2$. Let $(U,B,Q)$ be a non-degenerated quadratic space over $\mathbb{Q}_2$ and $W$ be a non-degenerated $\mathbb{Z}_2$-lattice on $U$. We say that $W$ has {\it even order} (resp. {\it odd order}) if $\ord_2(Q(\x))$ is even (resp. odd) for each primitive vector $\x\in W$ which gives rise to a $\mathbb{Z}_2$-integral symmetry of $W$.
       
       Suppose now that a non-degenerated lattice $M$ over $\mathbb{Z}_2$ has the following Jordan splitting
		$$M\cong 2^{r_1}M_{(1)}\perp 2^{r_2}M_{(2)}\perp\cdots\perp 2^{r_s}M_{(s)},$$
		where each $M_{(i)}$ is a unimodular lattice over $\mathbb{Z}_2$ and $0=r_1<r_2<\cdots<r_s$. In this paper, we will use the following result due to Earnest \cite[1.6]{Earnest78}.
		
		\begin{lemma}\label{Lem. spinor norm of Earnest}
			Let notations be as above. Then the following results hold. 
			
			{\rm (i)} Suppose that at least one Jordan component of $M$ has rank $\ge3$. Then 
			$$\theta_{\mathbb{Q}_2}\left(M\right)=
			\begin{cases}
				\mathbb{Z}_2^{\times}\mathbb{Q}_2^{\times2} & \mbox{if}\ 2^{r_k}M_{(k)} \ \text{has odd order for each}\ k=1,2,\cdots,s,\\
				\mathbb{Q}_2^{\times} & \mbox{otherwise.}
			\end{cases}$$
		
		    {\rm (ii)} If $\rank(M_{(k)})\le 2$ for each $k=1,2,\cdots,s$ and at least one component, say $M_{(k_0)}$ is binary, then $\theta_{\mathbb{Q}_2}(M)\neq \mathbb{Q}_2^{\times}$ if and only if one of the following three cases occur:
		    \begin{itemize}
		    	\item all Jordan components have odd order;
		    	
		    	\item all Jordan components have even order;
		    	
		    	\item whenever $\rank(M_{(k)})=2$, then 
		    	$$M_{(k)}\cong\left(\begin{array}{cc}
		    		\alpha & 1\\
		    		1      & 2\beta
		    	\end{array}\right)$$
	    	for some $\alpha,\beta\in\mathbb{Z}_2^{\times}$ and the associated quadratic space of $M_{(k)}$ is isotropic. Also, for any unary component, say $M_{(j)}\cong \langle c_j\rangle$ with $c_j\in\mathbb{Z}_2^{\times}$, the Hilbert symbol 
	    	$$\left(2^{r_j-r_{k_0}}\alpha c_j,-d(M_{(j)})\right)_2=1.$$
	    	Meanwhile, $r_{l+1}-r_l\ge 4$ for any $l=1,2,\cdots,s-1$.
		    \end{itemize}
		\end{lemma}

	We next briefly introduce some basic facts on spinor exceptional integers. Let $(V,B,Q)$ be a non-degenerated quadratic space over $\mathbb{Q}$ with $\dim_{\mathbb{Q}}(V)=3$, and let $L$ be a non-degenerated ternary lattice on $V$. An integer $t$ is called a {\it spinor exception} of $\gen(L)$ if $t\in Q(\gen(L))$ and $t$ is represented by exactly half of the spinor genera in $\gen(L)$. The arithmetic properties of spinor exception have been investigated extensively. Readers may refer to \cite{Duke,Earnest84,RSP80,RSP20}. In this paper, we need the following result due to R. Schulze-Pillot \cite{RSP20}.
	
	\begin{lemma}\label{Lem. RSP}
		Let $(V,B,Q)$ be a positive definite quadratic space over $\mathbb{Q}$ with $\dim_{\mathbb{Q}}(V)=3$, and let $L$ be a non-degenerated ternary lattice on $V$. Let $r\in\mathbb{Z}^+$ and let
		$$Q_r(\gen(L))=\left\{t\in\mathbb{Z}:\ t\in Q(\gen(L))\ \text{and}\ \ord_p(t)\le r\ \text{for any}\ p\in T\right\},$$
		where 
		$$T=\{p\in\Omega_{\mathbb{Q}}\setminus\{\infty\}:\ V_p\ \text{is anisotropic}\}.$$ 
		Then there is a constant $c$, depending on $r$ and the level of $L$, such that all $t\in  Q_r(\gen(L))$ with $t\ge c$ are represented by all lattices in $\gen(L)$ unless one of the following conditions are satisfied:
		
		{\rm (i)} $t$ is a spinor exception of $\gen(L)$;
		
		{\rm (ii)} $t/p^2$ is a spinor exception of $\gen(L)$ for some prime $p$ that is inert in the quadratic extension $\mathbb{Q}(\sqrt{-t\cdot d(L)})/\mathbb{Q}$. 
	\end{lemma}
	
	We conclude this section with the following result obtained by Sun \cite[Lemma 2.3]{Sun22}.
	
	\begin{lemma}\label{Lem. a lemma of sun to determine positivity}
		Let $a,b,c,d,m,n$ be non-negative real numbers with $a^2+b^2+c^2+d^2\neq 0$. Suppose that $x,y,z,w$ are real numbers satisfying 
		$$\begin{cases}
			x^2+y^2+z^2+w^2=n,\\
			ax+by+cz+dw=m.
		\end{cases}$$
	If
	$$m\ge \sqrt{\left((a^2+b^2+c^2+d^2)-\min\left(\left\{a^2,b^2,c^2,d^2\right\}\setminus\{0\}\right)\right)\cdot n},$$
	then all the numbers $ax,by,cz,dw$ are non-negative. 
	\end{lemma}

	\section{Proof of the main result}
	\setcounter{lemma}{0}
	\setcounter{theorem}{0}
	\setcounter{equation}{0}
	\setcounter{conjecture}{0}
	\setcounter{remark}{0}
	\setcounter{corollary}{0}

	{\noindent \bf Proof of Theorem \ref{Thm. A}}. Without loss of generality, we may assume $a\le b$. Let $(V,B,Q)$ be a positive definite quadratic space over $\mathbb{Q}$ with 
	$$V=\mathbb{Q}\e_1+\mathbb{Q}\e_2+\mathbb{Q}\e_3,$$
	where $\e_1,\e_2,\e_3$ is an orthogonal basis of $V$ with $Q(\e_1)=1$ and $Q(\e_2)=\Q(\e_3)=a^2+b^2$. Recall that $\gcd(a,b)=1$ and $a^2+b^2=p^r$ or $2p^k$, where $p$ is an odd prime, $r\in\mathbb{Z}^+$ and $k\in\mathbb{N}$. Moreover, one can easily verify that $p\equiv 1\pmod{4}$. 
	
	Let 
	$$L=\mathbb{Z}\e_1+\mathbb{Z}\e_2+\mathbb{Z}\e_3$$
	be a non-degenerated lattice on $V$, that is, 
	$$L\cong\langle1\rangle\perp\langle a^2+b^2 \rangle\perp \langle a^2+b^2 \rangle.$$
	For any integer $m$, define 
	$$l_{a,b,m}(n)=(a^2+b^2)n-m^2.$$
	
	We next consider the representations of $l_{a,b,m}(n)$ by the lattice $L_v$, where $v\in\Omega_{\mathbb{Q}}$. 
    
    {\bf (i)} For the archimedean place $v=\infty\in \Omega_{\mathbb{Q}}$ and each $n\in\mathbb{N}$, clearly
    \begin{equation}\label{Eq. representation of LR}
    	\left\{l_{a,b,m}(n):\ 0\le m\le \sqrt{(a^2+b^2)n}\right\}\subseteq Q(L_{\infty}).
    \end{equation}
	
	{\bf (ii)} For the place $v=q\in\Omega_{\mathbb{Q}}\setminus\{2,p,\infty\}$, the lattice $L_q$ is a unimodular lattice over $\mathbb{Z}_q$. For any $\varepsilon\in\mathbb{Z}_q^{\times}$, the extension $\mathbb{Q}_{q}(\sqrt{\varepsilon})/\mathbb{Q}_q$ is unramified. Thus, by the local class field theory,  
	$$\Norm_{\mathbb{Q}_q(\sqrt{\varepsilon})/\mathbb{Q}_q}(\mathbb{Q}_q(\sqrt{\varepsilon})^{\times})=
	\begin{cases}
		\mathbb{Q}_q^{\times} & \mbox{if}\ \varepsilon\in\mathbb{Z}_{q}^{\times2},\\
		\{q^{2\delta}\cdot u:\ \delta\in\mathbb{Z}\ \text{and}\ u\in\mathbb{Z}_q^{\times}\} & \mbox{if}\ \varepsilon\not\in\mathbb{Z}_{q}^{\times2}.
	\end{cases}$$
	By this, we obtain 
	$$S_q(V_q)=(-1,-1)_q=1,$$ 
	and hence $V$ is isotropic by Lemma \ref{Lem. isotropic}. Using Lemma \ref{Lem. unimodular lattices over Zp} again, 
	$$L_q\cong \left(\begin{array}{cc}
		0 & 1\\
		1 & 0
	\end{array}\right)\perp\langle-1\rangle.$$
	This implies $Q(L_q)=\mathbb{Z}_q$. Hence for any $n\in\mathbb{N}$, 
	\begin{equation}\label{Eq. representation of Lq}
		\left\{l_{a,b,m}(n):\ m\in\mathbb{Z}\right\}\subseteq Q(L_q).
	\end{equation}
	In addition, as $L_q$ is a ternary unimodular lattice over $\mathbb{Z}_q$, we have 
	\begin{equation}\label{Eq. spinor norm of Lq}
		\theta_{\mathbb{Q}_q}\left(L_q\right)=\mathbb{Z}_q^{\times}\mathbb{Q}_q^{\times2}.
	\end{equation}
	
	{\bf (iii)} For the place $v=p$, note that
	$$L_p\cong \langle 1 \rangle \perp \langle p^r\eta \rangle \perp \langle p^r\eta \rangle,$$
	where $r=\ord_p(a^2+b^2)\ge1$ and $\eta=(a^2+b^2)/p^r\in\mathbb{Z}_p^{\times}$. Since $p\equiv 1\pmod 4$, we have $-1\in\mathbb{Z}_p^{\times 2}$. By \cite[92:1]{OTO}, 
	$$\langle p^r\eta \rangle \perp \langle p^r\eta \rangle \cong \left(\begin{array}{cc}
		0   & p^r\\
		p^r & 0
	\end{array}\right).$$
	Thus, 
	$$L_p\cong  \langle 1 \rangle \perp \left(\begin{array}{cc}
		0   & p^r\\
		p^r & 0
	\end{array}\right).$$
	This implies that 
	$$Q(L_p)=\left\{x^2+p^ry:\ x,y\in\mathbb{Z}_p\right\}.$$
	By this and noting that $-1\in\mathbb{Z}_p^{\times 2}$, for any $n\in\mathbb{N}$, we obtain 
	\begin{equation}\label{Eq. representation of Lp}
		\left\{l_{a,b,m}(n):\ m\in\mathbb{Z}\right\}\subseteq Q(L_p).
	\end{equation}
	On the other hand, as $\langle\eta\rangle \perp \langle\eta\rangle$ is a binary unimodular lattice on $\mathbb{Z}_p$, we have 
	\begin{equation}\label{Eq. spinor norm of Lp}
		\mathbb{Z}_p^{\times}\mathbb{Q}_p^{\times2}=\theta_{\mathbb{Q}_p}\left(\langle\eta\rangle \perp \langle\eta\rangle\right)=\theta_{\mathbb{Q}_p}\left(\langle p^r\eta\rangle \perp \langle p^r\eta\rangle\right)\subseteq \theta_{\mathbb{Q}_p}\left(L_p\right).
	\end{equation}

	{\bf (iv)} We next consider the dyadic place $v=2$. Note that $\mathbb{Q}_2(\sqrt{-1})/\mathbb{Q}_2$ is a totally ramified extension. By the local class field theory, it is known that 
	\begin{equation}\label{Eq. local norm of Q2}
			\Norm_{\mathbb{Q}_2(\sqrt{-1})/\mathbb{Q}_2}(\mathbb{Q}_2(\sqrt{-1})^{\times})=
		\left\{2^{\delta}\cdot u:\ \delta\in\mathbb{Z}\ \text{and}\ u\in\mathbb{Z}_2^{\times}\ \text{with}\ u\equiv 1\pmod{4\mathbb{Z}_2}\right\}.
	\end{equation}

    We next consider the following two cases. 
    
    {\bf Case I}. $\ord_2(a^2+b^2)=0$. 
    
    In this case, 
    $$L_2\cong \langle 1 \rangle \perp \langle \mu \rangle \perp \langle \mu \rangle$$
    for some $\mu\in\mathbb{Z}_2^{\times}$ with $\mu\equiv 1\pmod{4\mathbb{Z}_2}$. By (\ref{Eq. local norm of Q2}), we see that 
    $$S_2(V_2)=1\neq (-1,-1)_2=-1,$$ 
    and hence $V_2$ is anisotropic by Lemma \ref{Lem. isotropic}. Applying Lemma \ref{Lem. unimodular lattices over Zp}, we obtain 
    $$L_2\cong\left(\begin{array}{cc}
    	2 & 1\\
    	1 & 2
    \end{array}\right)\perp \langle 3 \rangle\cong \langle1\rangle \perp \langle1\rangle \perp \langle1\rangle.$$
     This implies 
     \begin{equation}\label{Eq. values of Q(L2) in Case 1}
     	Q(L_2)=\mathbb{Z}_2\setminus\left\{4^{\delta}(8z+7):\ \delta\in\mathbb{N}\ \text{and}\ z\in\mathbb{Z}_2\right\}.
     \end{equation}
     
	By (\ref{Eq. values of Q(L2) in Case 1}) one can easily verify that following results hold.
	
	\begin{itemize}
		\item If $n\in\mathbb{N}$ with $n\equiv 1\pmod 4$, then 
		\begin{equation}\label{Eq. Case 1 representation of L2 when n=1 mod 4}
			\left\{l_{a,b,m}(n):\ m\in 2\mathbb{Z}\right\}\subseteq Q(L_2)\cap (4\mathbb{Z}+1)\neq\emptyset.
		\end{equation}
	    
	     \item If $n\in\mathbb{N}$ with $n\equiv 3\pmod 4$, then 
	     \begin{equation}\label{Eq. Case 1 representation of L2 when n=3 mod 4}
	     	\left\{l_{a,b,m}(n):\ m\in 2\mathbb{Z}+1\right\}\subseteq Q(L_2)\cap (4\mathbb{Z}+2)\neq\emptyset.
	     	\end{equation}
	     
	     \item Suppose $n\in\mathbb{N}$ with $\ord_2(n)=1$, i.e., $n\equiv 2\pmod 4$.  Then 
	     \begin{equation}\label{Eq. Case 1 representation of L2 when n=2 mod 4}
	     	\left\{l_{a,b,m}(n):\ m\in \mathbb{Z}\right\}\subseteq Q(L_2)\cap\left((4\mathbb{Z}+1)\cup (4\mathbb{Z}+2)\right)\neq\emptyset.
	     \end{equation}
	     
	     \item Suppose $n\in\mathbb{N}$ with $\ord_2(n)=2$, i.e., $n\equiv 4\pmod 8$. Then 
	     \begin{equation}\label{Eq. Case 1 representation of L2 when n=4 mod 8}
	     \left\{l_{a,b,m}(n):\ m\in 2\mathbb{Z}+1\right\}\subseteq Q(L_2)\cap (8\mathbb{Z}+3)\neq\emptyset.
	     \end{equation}
     
          \item Suppose $n\in\mathbb{N}$ with $\ord_2(n)=3$, i.e., $n\equiv 8\pmod {16}$. Then 
	     \begin{equation}\label{Eq. Case 1 representation of L2 when n=8 mod 16}
	     	\left\{l_{a,b,m}(n):\ m\in 4\mathbb{Z}\right\}\subseteq Q(L_2)\cap(16\mathbb{Z}+8)\neq\emptyset.
	     \end{equation}
	\end{itemize}

	Also, by Lemma \ref{Lem. spinor norm of Earnest}(i), 
	\begin{equation}\label{Eq. spinor norm of L2 in Case 1}
		\mathbb{Z}_2^{\times}\mathbb{Q}_2^{\times2}\subseteq\mathbb{Q}_2^{\times}=\theta_{\mathbb{Q}_2}\left(L_2\right).
	\end{equation}
	
	{\bf Case II}. $\ord_2(a^2+b^2)=1$.	
	
	Recall that $a^2+b^2=2p^k$ for some prime $p\equiv 1\pmod 4$ and some $k\in\mathbb{N}$ in this case. Then 
	$$L_2\cong\langle 1 \rangle \perp \langle 2\mu\rangle \perp \langle 2\mu\rangle $$
	for some $\mu\in\mathbb{Z}_2^{\times}$ with $\mu\equiv 1\pmod 4$. Noting that 
	$\e_1,\e_2+2\e_3, -2\e_2+\e_3$ is also a $\mathbb{Z}_2$-basis of $L_2$, we further obtain 
	$$L_2\cong \langle 1 \rangle \perp \langle 2\mu\rangle \perp \langle 2\mu\rangle \cong \langle 1 \rangle \perp \langle 10\mu\rangle \perp \langle 10\mu\rangle\cong \langle 1 \rangle \perp \langle 2\rangle \perp \langle 2\rangle.$$
	 By (\ref{Eq. local norm of Q2}), we see that 
	$$S_2(V_2)=1\neq (-1,-1)_2=-1,$$ 
	and hence $V_2$ is anisotropic by Lemma \ref{Lem. isotropic}. 
	
	Also, one can easily verify that the following results hold.
	\begin{itemize}
		\item If $n\in\mathbb{N}$ with $n\equiv 1\pmod 2$, then 
		\begin{equation}\label{Eq. Case 2 representation of L2 when n=1 mod 2}
			\left\{l_{a,b,m}(n):\ m\in 2\mathbb{Z}+1\right\}\subseteq Q(L_2)\cap (4\mathbb{Z}+1)\neq\emptyset.
		\end{equation}
	
	    \item Suppose $n\in\mathbb{N}$ with $\ord_2(n)=1$, i.e., $n\equiv 2\pmod 4$. Then 
	    \begin{equation}\label{Eq. Case 2 representation of L2 when n=2 mod 4}
	    	\left\{l_{a,b,m}(n):\ m\in 2\mathbb{Z}+1\right\}\subseteq Q(L_2)\cap (8\mathbb{Z}+3)\neq\emptyset.
	    \end{equation}
    
         \item Suppose $n\in\mathbb{N}$ with $\ord_2(n)=2$, i.e., $n\equiv 4\pmod 8$. Then 
         \begin{equation}\label{Eq. Case 2 representation of L2 when n=4 mod 8}
         	\left\{l_{a,b,m}(n):\ m\in 4\mathbb{Z}\right\}\subseteq Q(L_2)\cap (16\mathbb{Z}+8)\neq\emptyset.
         \end{equation}
         
         \item Suppose $n\in\mathbb{N}$ with $\ord_2(n)=3$, i.e., $n\equiv 8\pmod {16}$. Then
         \begin{align}\label{Eq. Case 2 representation of L2 when n=8 mod 16}
         	   &\left\{l_{a,b,m}(n):\ m\in 4\mathbb{Z}\ \text{and}\ \frac{m}{4}\equiv \frac{(a^2+b^2)n-16}{32}\pmod 2 \right\}\notag\\
               &\subseteq Q(L_2)\cap \left((64\mathbb{Z}+16)\cup(64\mathbb{Z}+32)\right)\neq\emptyset.
         \end{align}
	\end{itemize}

	In addition, by Lemma \ref{Lem. spinor norm of Earnest}(ii), 
	\begin{equation}\label{Eq. spinor norm of L2 in Case 2}
		\mathbb{Z}_2^{\times}\mathbb{Q}_2^{\times2}\subseteq\mathbb{Q}_2^{\times}=\theta_{\mathbb{Q}_2}\left(L_2\right).
	\end{equation}

	After the discussions of the above two cases, combining (\ref{Eq. spinor norm of L2 in Case 1}) and (\ref{Eq. spinor norm of L2 in Case 2}) with (\ref{Eq. spinor norm of Lq}) and (\ref{Eq. spinor norm of Lp}), by Lemma \ref{Lem. number of spinor genera} we obtain that there is a unique spinor genus in $\gen(L)$, i.e., 
		$$\spn(L)=\gen(L),$$
	and 
	\begin{equation}\label{Eq. the only anisotropic prime is 2}
		\left\{v\in\Omega_{\mathbb{Q}}\setminus\{\infty\}: \text{the associated quadratic space of}\ L_v\ \text{is anisotropic}\right\}=\left\{2\right\}.
	\end{equation}

	Now for any sufficiently large integer $n$ with $\ord_2(n)\le C$, we write $n=16^{\delta}\cdot n_1$, where 
	$$\delta=\max\left\{s\in\mathbb{N}:\ n\equiv 0\pmod{16^s}\right\}.$$
	As $\ord_2(n)\le C$, we see that $n_1$ is also a sufficiently large integer. Thus, we may assume that 
	$$\left((a^2+b^2)n_1\right)^{1/4}-\left(b^2n_1\right)^{1/4}>8.$$
	This implies that there is an integer $\alpha$ such that 
	\begin{equation}\label{Eq. there are 2 squares in the interval}
		(\alpha+i)^2\in\left[\left(b^2n_1\right)^{1/2}, \left((a^2+b^2)n_1\right)^{1/2}\right]
	\end{equation}
    for $i=0,1,2,\cdots,7$, that is, there are at least eight consecutive squares in the above interval. Using (\ref{Eq. there are 2 squares in the interval}) and by (\ref{Eq. representation of LR}), (\ref{Eq. representation of Lq}), (\ref{Eq. representation of Lp}), (\ref{Eq. Case 1 representation of L2 when n=1 mod 4})--(\ref{Eq. Case 1 representation of L2 when n=8 mod 16}) and (\ref{Eq. Case 2 representation of L2 when n=1 mod 2})--(\ref{Eq. Case 2 representation of L2 when n=8 mod 16}), we can choose a positive square $m^2\in [(b^2n_1)^{1/2}, ((a^2+b^2)n_1)^{1/2}]$ satisfying the following three conditions:  
    \begin{itemize}
    	\item $l_{a,b,m^2}(n_1)\in Q(\gen(L))$,
    	\item $\ord_2(l_{a,b,m^2}(n_1))\le 5$,
    	\item $\gcd(m,p)=1$.
    \end{itemize}
	This, together with $\spn(L)=\gen(L)$, (\ref{Eq. the only anisotropic prime is 2}) and Lemma \ref{Lem. RSP}, implies that 
	$$l_{a,b,m^2}(n_1)\in Q(L).$$
	Thus, there exist $u\in\mathbb{Z}$ and $z,w\in\mathbb{N}$ such that 
	\begin{equation}\label{Eq. representation involving uzw}
		(a^2+b^2)n_1-m^4=u^2+(a^2+b^2)z^2+(a^2+b^2)w^2.
	\end{equation}

    As $\gcd(a,b)=1$, there exist $s,t\in\mathbb{Z}$ such that 
    $$m^2=as+bt.$$ 
    Also, by this we obtain 
	\begin{equation}\label{Eq. congruence involving s and t}
		\left(\pm(at-bs)\right)^2\equiv -(as+bt)^2\equiv -m^4\pmod{a^2+b^2}.
	\end{equation}

We first consider the case $a^2+b^2>2$. Recall that the multiplicative group
	$$\left(\mathbb{Z}/(a^2+b^2)\mathbb{Z}\right)^{\times}=\left\{c \mod{(a^2+b^2)\mathbb{Z}}:\ \gcd(c,a^2+b^2)=1\right\}$$
	is a cyclic group of even order. 
	
	Suppose $\gcd(m,a^2+b^2)=1$. Then by the property of even cyclic group and (\ref{Eq. congruence involving s and t}), 
	$$u\equiv \pm(at-bs)\pmod{a^2+b^2}$$
	are exactly all the solutions of the equation 
	$$ u^2\equiv -m^4\pmod{a^2+b^2}$$
	in the even cyclic group $(\mathbb{Z}/(a^2+b^2)\mathbb{Z})^{\times}$. By this,  (\ref{Eq. representation involving uzw}) and replacing $u$ by $-u$ if necessary, we may set 
	$$u=(a^2+b^2)u'+(at-bs)$$
	for some $u'\in\mathbb{Z}$. Using this and (\ref{Eq. representation involving uzw}) again, we deduce that 
	$$n_1=(s-bu')^2+(t+au')^2+z^2+w^2.$$
	Letting $x=s-bu'$, $y=t+au'$ and noting that $m^2=as+bt$, we obtain 
	$$\begin{cases}
		x^2+y^2+z^2+w^2=n_1,\\
		ax+by=m^2.
	\end{cases}$$
	
	Suppose $\gcd(m,a^2+b^2)>1$. Then by the choice of $m$, in this case, we have $2\nmid ab$, $\ord_2(a^2+b^2)=1$ and $\gcd(m,a^2+b^2)=2$. As $a^2+b^2>2$, as mentioned before, the multiplicative group 
	$$\left(\mathbb{Z}/p^k\mathbb{Z}\right)^{\times}=\left\{c \mod{p^k\mathbb{Z}}:\ \gcd(c,p^k)=1\right\}$$
	is also a cyclic group of even order, where 
	$$p^k=(a^2+b^2)/2.$$
	As $\gcd(m,p)=1$, we see that 
	$$u\equiv \pm(at-bs)\pmod{p^k}$$
	are exactly all the solutions of the equation 
	$$ u^2\equiv -m^4\pmod{p^k}$$
	in the even cyclic group $(\mathbb{Z}/p^k\mathbb{Z})^{\times}$. Thus, by (\ref{Eq. representation involving uzw}) and replacing $u$ by $-u$ if necessary, we may assume 
	$$u=\frac{a^2+b^2}{2}u'+(at-bs)$$
	for some $u'\in\mathbb{Z}$. Note that $at-bs$ and $u$ are both even and $(a^2+b^2)/2$ is odd. We therefore obtain that $u'=2u''$ is an even integer, that is,  
	$$u=(a^2+b^2)u''+(at-bs).$$
	This, together with (\ref{Eq. representation involving uzw}), implies that 
	$$n_1=(s-bu'')^2+(t+au'')^2+z^2+w^2.$$
	Letting $x=s-bu''$ and $y=t+au''$ and noting that $as+bt=m^2$ again, we also have 
	$$\begin{cases}
		x^2+y^2+z^2+w^2=n_1,\\
		ax+by=m^2.
	\end{cases}$$
	
	Now we consider the case $a=b=1$. As $s+t=m^2$, (\ref{Eq. representation involving uzw}) implies 
	$$u\equiv m^2\equiv s+t\equiv s-t\pmod 2.$$
	Thus, by replacing $u$ by $-u$ if necessary, we may set 
	$$u=2u'''+(t-s)$$
	for some $u'''\in\mathbb{Z}$. 
	By this and  (\ref{Eq. representation involving uzw}) again, we obtain 
	$$2n_1-m^4(2u'''+(t-s))^2+2z^2+2w^2,$$
	which implies that 
	$$n_1=(s-u''')^2+(t+u''')^2+z^2+w^2.$$
	Setting $x=s-u''', y=t+u'''$, we obtain 
	$$\begin{cases}
		x^2+y^2+z^2+w^2=n_1,\\
		x+y=m^2.
	\end{cases}$$
	
	By the above discussions, for $n_1$, we can always find a square $m^2\in [(b^2n_1)^{1/2}, ((a^2+b^2)n_1)^{1/2}]$, $x,y\in\mathbb{Z}$ and $z,w\in\mathbb{N}$ such that 
	$$\begin{cases}
		x^2+y^2+z^2+w^2=n_1,\\
		ax+by=m^2.
	\end{cases}$$
    Now applying Lemma \ref{Lem. a lemma of sun to determine positivity}, we immediately obtain that $x,y\in\mathbb{N}$. Recall that $n=16^{\delta}\cdot n_1$. We thus directly obtain 
    $$\begin{cases}
    	(4^{\delta}x)^2+(4^{\delta}y)^2+(4^{\delta}z)^2+(4^{\delta}w)^2=n,\\
    	a(4^{\delta}x)+b(4^{\delta}y)=(2^{\delta}m)^2\in\mathcal{S}.
    \end{cases}$$
	
	In view of the above, we have completed the proof.\qed

	\Ack\ This research was supported by the Natural Science Foundation of China (Grant No. 12101321). The first author was also supported by the Natural Science Foundation of the Higher Education Institutions of Jiangsu Province (Grant No. 25KJB110010).

    \medskip\noindent {\bf Declaration of competing interest}
    
    The authors declare that they have no conflict of interest.
	
	\medskip\noindent {\bf Data avilability}

    No data was used for the research described in the article.

\end{document}